%% file: main.tex
\documentclass[draft]{agujournal2019}

\usepackage{amssymb, amsmath, amscd}
\usepackage{graphicx}
\usepackage{xfrac} 

\usepackage{url} 
\draftfalse

\journalname{Journal of Advances in Modeling Earth Systems (JAMES)}

\begin{document}

\title{An energy consistent discretization of the nonhydrostatic equations in primitive variables}

%
%



\authors{
Mark A. Taylor\affil{1},
Oksana Guba\affil{1},
Andrew Steyer\affil{1},
Paul A. Ullrich\affil{2}, 
David M. Hall\affil{3},
Christopher Eldred\affil{4}
}

\affiliation{1}{Computational Science, Sandia National Laboratories, Albuquerque,
New Mexico, USA.}
\affiliation{2}{Department of Land, Air and Water Resources, 
University of California, Davis, California, USA.}
\affiliation{3}{NVIDIA, Santa Clara, California, USA}
\affiliation{4}{Univ. Grenoble Alpes, Inria, CNRS, Grenoble INP, LJK, 38000 Grenoble, France}




\correspondingauthor{Mark Taylor}{mataylo@sandia.gov}



\begin{keypoints}
\item	
We give a discrete Hamiltonian formulation of the nonhydrostatic 
equations in primitive variables.
\item	
The formulation supports mass or height terrain following coordinates.
\item	
The Lorenz staggered vertical discretization obeys a derivative product rule.
\end{keypoints}

%
%


\begin{abstract}
We derive a formulation of the nonhydrostatic equations in 
spherical geometry with a Lorenz staggered vertical discretization.  The combination conserves a 
discrete
energy in exact time integration when coupled with a mimetic horizontal
discretization.  The formulation is a version of \citeA{DubosTort2014} rewritten in terms of 
primitive variables.  It is valid for terrain following mass or
height coordinates and for both Eulerian or vertically Lagrangian discretizations.
The discretization relies on an extension to \citeA{SB81} vertical differencing
which we show obeys a discrete derivative product rule.  
This product rule allows us to simplify the treatment of the 
vertical transport terms.   Energy conservation is obtained via a term-by-term balance in the kinetic,
internal and potential energy budgets, ensuring an energy-consistent
discretization with no spurious sources of energy.  We demonstrate convergence with respect to
time truncation error in a spectral element code with a HEVI IMEX
timestepping algorithm.  
\end{abstract}

\newcommand\mpar[1]{
\marginpar{\renewcommand{\baselinestretch}{1}\raggedleft
\tt\tiny\hspace{0pt}#1}%
}

\newcommand\PAU{\textcolor{red}}

%
%

\section*{Plain Language Summary}
Energy consistent discretizations have proven useful in guiding the development
of numerical methods for simulating fluid dynamics.  They ensure that the
discrete method does not have any spurious sources of energy which can lead
to unstable and unrealistic simulations.   Here we provide an energy consistent
discretization of the equations used by global models of 
the Earth's atmosphere.  The discretization is written in terms of 
standard variables in spherical coordinates and supports a wide variety of 
terrain following vertical
coordinates.  It can be used with any horizontal discretization that has
a discrete version of the integration-by-parts identity.

\input{introduction}

\input{narrative}

%
%
%
%
%
%
%
%

\acknowledgments

HOMME-NH is part of the open source E3SM project with source code
publically available from
\url{https://github.com/E3SM-Project/E3SM}. Instructions for
configuring and compiling standalone HOMME-NH are in the README Files in
\path{components/homme}.  Once the model has been configured, the
\path{checkE.sh} script in
\path{dcmip_tests/dcmip2016_test1_baroclinic_wave/theta-l} will
generate all necessary input data and run the simulations described in
this work.

This research was supported as part of the Energy Exascale Earth System Model (E3SM) project, funded by the U.S. Department of Energy, Office of Science, Office of Biological and Environmental Research.  
Sandia National Laboratories is a multimission laboratory managed and
operated by National Technology \& Engineering Solutions of Sandia,
LLC, a wholly owned subsidiary of Honeywell International Inc., for
the U.S. Department of Energy’s National Nuclear Security
Administration under contract DE-NA0003525.
This paper describes objective technical results and analysis. Any
subjective views or opinions that might be expressed in the paper do
not necessarily represent the views of the U.S. Department of Energy
or the United States Government.

Christopher Eldred was supported by the French National Research Agency through contract ANR-14-CE23-0010 (HEAT).


%
%

\bibliography{main}

%
%
%
%
%

\end{document}

%% file: introduction.tex
\section{Introduction}

We present a discrete formulation of the nonhydrostatic equations
designed for use in a global model of the Earth's atmosphere.  It has
been implemented in the version of the High Order Method Modeling
Environment (HOMME), \cite{Dennis05,Dennis12} used by the Energy
Exascale Earth System Model (E3SM), \cite{Golaz2019}.  We follow a
Hamiltonian approach so that the formulation will be energetically
conservative and consistent in the sense of \citeA{GassmannHerzog2008}
when coupled with a mimetic horizontal and vertical discretizations.
Energy consistency is obtained via a term-by-term balance in the
discrete kinetic, internal and potential energy budgets, ensuring an
energy conserving discretization with no spurious energy sources or
sinks.




Recent examples of energy consistent global atmosphere models running
on unstructured grids include \cite{Taylor11,Gassmann2013,Dubos2015}.
Here we extend the discretization from \citeA{Taylor11} to the
nonhydrostatic equations using a formulation modeled after
\citeA{DubosTort2014}. The main difference is that we formulate the
Hamiltonian in terms of the commonly used primitive variables in
spherical coordinates.  This may simplify implementation in existing
modeling frameworks.  The discretization supports terrain following
mass or height coordinates \cite{Kasahara1974}, which differ
only in their boundary conditions and one diagnostic equation.
It retains a 2D vector
invariant form for horizontal transport and advective form for
vertical transport, facilitating both Eulerian and
vertically Lagrangian \cite{sjlin04} coordinates in a single code
base.  In the vertically Lagrangian (mass or height) coordinate, the
vertical coordinate surfaces float in such a way as to ensure no mass
flux through the surfaces.
In the Eulerian (mass or height) coordinate, the transformed equations
are solved in an Eulerian reference frame.

For our vertical discretization, we use \citeA{SB81} vertical
differencing with a \citeA{Lorenz1960} staggering.  We introduce an
extension to SB81 which obeys a discrete derivative product rule
in addition to the well known SB81 discrete integral properties.  This
product rule ensures the discrete vertical transport terms remain
energetically neutral despite being in non-Hamiltonian form.  The
discretization is energetically consistent when used with the shallow
atmosphere approximation, total air mass density as a prognostic
variable and general moist equations of state, but with a common
approximation in the equation for virtual potential temperature.

The formulation is appropriate for any mimetic spatial
discretization in
spherical geometry such as
\citeA{thuburn2009numerical,ringler2010unified,TayFour10,cotter2012mixed,thuburn2015primal,
Gassmann2018}.  Mimetic discretizations are those that mimic key
integration and annihilator properties of the divergence, gradient and curl
operators 
\cite{Samarskii81,Nicolaides92, ShashkovSteinberg95,bh2006}. 
  They are commonly used in 
   order to preserve the Hamiltonian structure of the continuum
equations \cite{Salmon04,Salmon07,GassmannHerzog2008}. 
Here we use a $p=3$ collocated spectral element discretization due to its
retention of 4th order accuracy on unstructured spherical grids and its
explicit mass matrix which allows for the use of efficient 
horizontal explicit vertically implicit (HEVI) timestepping.
This efficiency comes at the cost of erratic behavior of grid scale waves
\cite{melvin2012,ainsworth2014}
which are damped
as in \citeA{ullrich2018}. 
We demonstrate convergence with respect to time truncation error with
a 5 stage Runge-Kutta IMEX algorithm.  

%% file: narrative.tex
\newcommand{\llangle}{\left\langle}
\newcommand{\rrangle}{\right\rangle}
\renewcommand{\vec}{\mathbf}
\newcommand{\grad}{\nabla}
\renewcommand{\div}{\nabla\cdot}
\newcommand{\curl}{\nabla\times}
\newcommand{\vvec}{\vec{v}}   
\newcommand{\uvec}{\vec{u}}   
\newcommand{\khat}{\hat k}
\newcommand{\dpids}{\frac{\partial \pi}{\partial s}}
\newcommand{\tdpids}{\tfrac{\partial \pi}{\partial s}}

\newcommand{\DS}{{\Delta s}}
\newcommand{\half}{{\sfrac12}}  

\section{Continuum equations in terrain following coordinates}
\label{S:transform}

Following \citeA{Kasahara1974} and \citeA{Laprise92} we write
the equations of motion in a terrain following coordinate and
make use of a pseudo-density based on the hydrostatic pressure.
Our physical domain will be Earth's 3D atmosphere in spherical coordinates with
latitude $\varphi$, longitude $\lambda$, and radial coordinate $r$.  The altitude from an arbitrary reference level $r_0$ is denoted $z = r - r_0$.  A terrain following coordinate $s$ is defined so that the bottom boundary corresponds to a surface of constant $s$.
We denote
the velocity $\vvec = (u,v,w)$  with $\uvec = (u,v)$ denoting the
the horizontal velocity tangent to constant $z$ surfaces.
For vertical velocity, we use both 
$\dot z = w = Dz/Dt$ and $\dot s = Ds/Dt$, with $D/Dt$ the
standard material derivative.
We take $\rho$ as the total air density,
$\khat$ the radial unit vector, and $\phi=gz$ as
the geopotential with constant gravitational acceleration $g$.  

We assume $z$ is a monotone function of $s$, so
the pseudo-density $\partial \pi/\partial s$ can be defined in terms
of the hydrostatic pressure $\pi$ and density,
 \[
 \dpids\frac{\partial s}{\partial z} = \frac{\partial \pi}{\partial z} = -\rho g
\]
with
\[
\rho  = - \dpids \Big/  \frac{\partial \phi}{\partial s}.
\]
Vertical integrals of mass and other density weighted quantities
are transformed to terrain following coordiantes via
\[
\int_{z_\text{surf}}^{z_\text{top}} \rho X   \,dz  =
\frac1g \int_{s_\text{top}}^{s_\text{surf}} \dpids  X  \,ds.
 \]
We denote the three-dimensional physical coordinate-independent differential operators 
in the usual notation,
$
\div \vvec, \grad p, \curl \vvec
$
and use 
\[
\grad_s p,
\quad
\nabla_s \times w \hat k,
\quad
\nabla_s \cdot \uvec,
\quad \nabla_s \times  \uvec
\]
to denote the two-dimensional operators on $s$ surfaces
\cite{ullrich2017dcmip}.

\subsection{Pressure gradient formulation}
\label{S:pgrad_formulations}

We can separate the pressure gradient term into vertical and horizontal components,
written in terrain following coordinates as
\begin{align}
 \frac1\rho  \grad p \cdot \khat   &= 
g \frac1\rho  \frac{\partial p}{\partial s} \Big/  \frac{\partial \phi}{\partial s}
=
-g   \mu
 \\
 \frac1\rho  \{ \grad p \}_\text{h}  
&=
\frac1\rho  \grad_s p + \mu \grad_s \phi 
\end{align}
where for compactness we introduce 
\[
\mu =   \frac{\partial p}{\partial s}  \Big/ \dpids.
\]
We use the usual definition of 
virtual Exner pressure
$\Pi = (p/p_0)^{R/c_p}$ and virtual potential temperature
$\theta_v = p / (\rho R \Pi)$.  These are defined
in terms of $p$ and $\rho$ from the moist atmosphere but
using the ideal gas constant for dry air, $R$, and the 
specific heat at constant pressure for dry air, $c_p$.
These definitions do not imply any assumptions
regarding the equation of state and are used to write
the pressure gradient as 
\[
\frac1\rho  \grad_s p = c_p \theta_v \grad_s \Pi.
\]

\subsection{Shallow atmosphere equations}
We start with a standard form of the nonhydrostatic equations with the shallow atmosphere approximation in terrain following coordinates.  
Our goal is to derive a energetically consistent
discretization for the adiabatic terms, so we initially 
neglect dissipation, sources and sinks
of moisture and other forcing terms.  In the full model, dissipation will be added through choice of timestepping and vertical remap algorithms and horizontal hyperviscosity.    
Following \citeA{Laprise92}, we write the equations as
\begin{align}
  \frac{\partial \uvec}{\partial t} + (\nabla_s \times  \uvec + f \khat )  \times \uvec 
+\frac12 \grad_s \uvec^2 
+ \dot s \frac{\partial \uvec}{\partial s} 
+c_p \theta_v  \grad_s \Pi + \mu \grad_s \phi
&= 0
\label{E:adiabatic1}
\\
  \frac{\partial w}{\partial t} +
\uvec \cdot \grad_s w + \dot s \frac{\partial w}{\partial s}
+g 
- g \mu
&= 0
\label{E:adiabatic2}
\\
  \frac{\partial \phi}{\partial t} +
  \uvec \cdot \grad_s \phi + \dot s \frac{\partial \phi}{\partial s} - gw &= 0
  \label{E:adiabatic3}
\\
\dfrac{\partial \theta }{\partial t}   + \uvec \cdot \nabla_s  \theta   + 
\dot s \dfrac{\partial  \theta}{\partial s}      &= 0
  \label{E:adiabatic4}
\\
\dfrac{\partial }{\partial t} \left( \dpids \right)  + \nabla_s \cdot  \left( \dpids  \uvec \right)    + 
\dfrac{\partial }{\partial s} \left(  \dpids  \dot s \right) &= 0
  \label{E:adiabatic5}
\end{align}
with Coriolis term $f\khat$.  

\subsection{Equation of state}
To close the equations we need an equation of state in order to
compute $p$, which is needed to derive $\Pi$ and $\theta_v$.  
As written above, we can
support a very general equation of state
via the Gibbs potential approach recommended in
\citeA{Thuburn2017}.  Given a Gibbs potential $G(p,T,q_i)$, with
temperature $T$ and $q_i$ representing various forms of water, 
the algebraic system
\begin{align}
    G_T(p_0,\theta,q_i) &= G_T(p,T,q_i) \\
    \frac1\rho &= G_p(p,T,q_i)
\end{align}
then defines $p$ and $T$ implicitly as functions of $\rho,\theta$ and $q_i$.
This system can be solved efficiently with Newton iteration.

Here we simplify the equations substantially by assuming that
\[
\frac{D \theta_v}{D t} = 0.
\]
This is a common assumption for moist air.  It 
neglects small terms, the magnitude of which 
are discussed in \citeA{staniforth2006}.
Under this approximation, we can prognose $\theta_v$ instead of $\theta$
and then the relation $\theta_v = p/(\rho R \Pi) $ is sufficient to close
the system.  The water species $q_i$ then become passive tracers during
the adiabatic part of the dynamics represented by
\eqref{E:adiabatic1}-\eqref{E:adiabatic5}.
The equation of state is only needed when coupling with other
aspects of the model, such as to compute $T$ needed by the physical 
parameterizations.
The $D\theta_v/Dt=0$ approximation can be made independently of the
equation of state, but it is more consistent to introduce an
approximate Gibbs potential for which $D\theta_v/Dt=0$ and
use this potential in all equation of state calculations
\cite{Thuburn2017}.

For the simulations presented here, we use an equation of
state for dry air, $\rho_d$, and water vapor represented by
specific humidity $q$.
We assume $q$ is prognosed by the model's transport scheme.
The Gibbs potential for this equation of state,
e.g. \citeA[Appendix A]{Vallis2017}, leads to 
$ p = \rho R^* T $
with
$\rho_d = \rho(1-q)$, $R^*  = R  + \left( R_v - R  \right) q  $,
$c^*_p = c_p +  \left( c_{pv} - c_p  \right) q $
where $R$ and $c_p$ are the dry air constants defined above,
$R_v$ is the ideal gas constants for water vapor and
$c_{pv}$ is the specific heat at constant pressure for water vapor.
For this equation of state, we have
\[
T  = \theta  \left(\frac{p}{p_0}\right)^{R^*/c_p^*}  
\]
and $D\theta/Dt = 0$ in the absence of sources or sinks of moisture ($Dq/Dt = 0$).

We can introduce an approximate Gibbs potential that is consistent
with $D\theta_v / D t = 0$ by taking 
$c_{pv} = (R_v/R) c_p $.  This representing a 14\% error in $c_{pv}$.
With this $c_{pv}$ approximation,  $R^*/c_p^* = R/c_p$, so that 
$R \theta_v = R^* \theta$ and then  $D\theta_v / D t = 0$ since
$DR^*/Dt = (R_v -R) Dq/Dt = 0$.

%
%

 \subsection{Quasi-Hamiltonian Form}
 We now give the final form of the continuum equations that
 will be discretized in \S~\ref{S:vert_disc}.  
We write the thermodynamic equation in conservation form instead
of advective form.  We also 
change the formulation slightly so that it is in
quasi-Hamiltonian form (See \ref{S:appendix}) which requires the gradient of the full 
kinetic energy in the momentum equation.
The result is
\begin{align}
  \frac{\partial \uvec}{\partial t} + (\nabla_s \times  \uvec + f\khat)  \times \uvec 
+\frac12 \grad_s \left( \uvec^2  + w^2 \right)
- w \grad_s w
+ \dot s \frac{\partial \uvec}{\partial s} 
+c_p \theta_v  \grad_s \Pi + \mu \grad_s \phi
&= 0
\label{E:uequation}
\\
  \frac{\partial w}{\partial t} +
\uvec \cdot \grad_s w + \dot s \frac{\partial w}{\partial s}
+g 
- g \mu
&= 0
\label{E:wequation}
\\
  \frac{\partial \phi}{\partial t} +
  \uvec \cdot \grad_s \phi + \dot s \frac{\partial \phi}{\partial s} - gw &= 0
  \label{E:phi}
\\
\dfrac{\partial \Theta }{\partial t}   + \grad_s \cdot (\uvec \Theta )
+    \dfrac{\partial}{\partial s} \left( \dot s \Theta \right)      &= 0
\\
\dfrac{\partial }{\partial t} \left( \dpids \right)  + \nabla_s \cdot  \left( \dpids  \uvec \right)    + 
\dfrac{\partial }{\partial s} \left(  \dpids  \dot s \right) &= 0
\label{E:continuity}
\end{align}
with
\begin{equation}
\Theta = \dpids \theta_v \qquad
\frac{\partial \phi}{\partial s}   = - R \Theta \frac{\Pi}{p}
\label{E:EOS}
\end{equation}
where in \eqref{E:EOS} we have rewritten $\theta_v=p/(\rho R \Pi) $ in terms of our prognostic
variables.

\subsection{Height coordinate}

In the height coordinate, our formulation is valid for any
general terrain following coordinate as long as $z$ is a monotone
function of $s$ and independent of time, such as \citeA{GalChen1975}
and generalizations
that allow for a quicker transition from 
terrain following to pure height such as  \citeA{Schar2002, Leuenberger2010, klemp2011}.

For boundary conditions at the model top and surface, we
use a free slip condition ($\dot s=0$).  Since 
$\partial \phi / \partial t = 0$ throughout the domain,
\eqref{E:phi} remains a valid equation if we take diagnostic
equation for  $\dot s$ as
\begin{equation}
\dot s  \frac{\partial \phi}{\partial s}  = g w - \uvec \cdot \grad_s \phi.
\label{E:sdot-height}
\end{equation}
Evaluating this equation at the model top and surface we derive
the $w$ boundary condition
\[
w  =  \uvec \cdot \frac1g \grad_s \phi
\]
which is combined with \eqref{E:wequation} to derive boundary conditions
for $p$ as well.

\subsection{Mass coordinate}
For the mass coordinate we follow \citeA{ArakawaLamb1977,SB81} and
define $s$ via 
\[
\pi = A(s) p_0 + B(s) p_s
\]
with $A(s)$ zero near the surface and
$B(s)$ zero near the model top.  
In this coordinate system, the diagnostic equation for $\dot s$ is
obtained by integrating \eqref{E:continuity} in the vertical and 
making use of $\partial \pi / \partial t =
B(s) \partial   \pi_s / \partial t$ to obtain
\begin{equation}
\dot s \dpids = 
 B(s) \int_{s_\text{top}}^1 \grad_s \cdot  \dpids \uvec \,ds  
 - \int_{s_\text{top}}^s \grad_s \cdot  \dpids \uvec \,ds
 \label{E:sdot-mass}
\end{equation}
For boundary conditions at the model top and surface, we
use a free slip condition ($\dot s=0$).
Since $\partial \phi / \partial t = 0$ at the model surface,
from \eqref{E:phi} we have a boundary condition on $w$,
\[
w  =  \uvec \cdot \frac1g \grad_s \phi
\]
which is combined with \eqref{E:wequation} to derive boundary conditions
for $p$ at the surface.   At the top of the model we have a free surface
with a constant pressure $p = A(s_\text{top}) p_0$.

Following \citeA{Laprise92}, the standard mass coordinate hydrostatic
approximation can be obtained by taking $p=\pi$ and replacing 
\eqref{E:phi} with the diagnostic equation for $\phi$ obtained by integrating \eqref{E:EOS}.

\subsection{Vertically Lagrangian coordinate}
We also consider the vertically Lagrangian coordinate from \citeA{sjlin04}.
We assume the vertical levels move so that there is no
flow through the levels ($\dot s = 0$) and the
vertical advection terms are dropped from the equations.
In both mass and height options, 
the diagnostic equation for $\dot s$ is no longer valid nor
is it needed.  Instead, \eqref{E:phi}
must be used as a prognostic equation for the level positions.

In practice, the floating vertical levels will eventually become too
unevenly spaced and every few timesteps we must remap the solution
onto a set of more uniform levels.  Here we consider energy
conservation only during the floating Lagrangian phase.  To remap, we
typically use a conservative monotone algorithm which will introduce
vertical dissipation.  To retain energy conservation in the vertically
Lagrangian case an energy conserving remap algorithm would have to be
adopted.

\section{Energy conservation}
\label{S:EnergyConservation}

We first show energy conservation and derive the energy budget
in the continuum.  We look at the budget term by term to
check that we have formulated the equations so that most
terms cancel via integration by parts, so that they will
also cancel when discretized
with a numerical method that has a discrete integration by parts
property.  As such we must be careful to only perform algebraic
steps which will also hold in the discrete system.  
In this derivation we note a few key terms that require additional
numerical properties in order to obtain energy conservation.
The analysis in the discrete case will
be identical but with integrals replaced by sums and additional
care must be taken to account for the vertical staggering.

The kinetic, internal and potential energy of our equations
in terrain following coordinates with integration with respect
to $ds$ is given by
\begin{equation}
K = \frac12 \dpids \vvec^2
\qquad
I = c_p  \Theta  \Pi +  \frac{\partial \phi}{\partial s} p  
+ \hat p  \phi_\text{top}
\qquad
P = \dpids \phi
\label{E:energy_pcoords}
\end{equation}
with $\hat p = p_\text{top}$ in the mass coodrinate and $\hat p = 0$ in the height
coordinate. In the mass coordiante, $p_\text{top}$
appears in the energy due to the free surface constant pressure boundary condition.

\subsection{Potential energy}
For the evolution of potential energy in the mass coordinate,
we multiply the continuity equation by $\phi$
and the $\phi$ equation by $\dpids$ and sum,
\begin{equation}
\frac{\partial }{\partial t}P =
-  \phi \nabla_s \cdot  \left( \dpids  \uvec \right)     
- \dpids \uvec \cdot \grad_s \phi
-   \phi \dfrac{\partial }{\partial s} \left(  \dpids  \dot s \right) 
- \dpids \dot s \frac{\partial \phi}{\partial s} 
+ \dpids gw 
\end{equation}
The first four terms will integrate to zero, requiring only
integration by parts, and the remaining term is the transfer of
potential energy to kinetic energy.  In the Lorenz staggered case,
some care will be needed in the averaging to ensure the first four
terms cancel.
This equation can be simplified in the Eulerian height coordinate
since $\partial \phi / \partial t=0$, but since this term is not zero
in the mass coordinate or the vertically Lagrangian height coordinate,
we retain this general form.

\subsection{Internal energy}
For internal energy we have
\[
 \dfrac{\partial }{\partial t} I = 
c_p  \Pi \dfrac{\partial \Theta}{\partial t} 
+ c_p  \Theta  \dfrac{\partial \Pi }{\partial t}
+ \dfrac{\partial }{\partial t} \left(   \frac{\partial \phi}{\partial s} p     \right)
+ \hat p \dfrac{\partial }{\partial t} (  \phi_\text{top}).  
\]
Assuming exact time integration, this can be simplified to
\[
  \dfrac{\partial }{\partial t} I 
= 
c_p  \Pi \dfrac{\partial \Theta}{\partial t}
+ p \dfrac{\partial }{\partial s}  \frac{\partial \phi}{\partial t} 
+ \hat p \dfrac{\partial }{\partial t} (  \phi_\text{top})
\]
We can eliminate the mass coordinate boundary term if we integrate the second term by
parts in the vertical, leading to
\[
  \dfrac{\partial }{\partial t} I 
= 
c_p  \Pi \dfrac{\partial \Theta}{\partial t}
- \dfrac{\partial p}{\partial s}  \frac{\partial \phi}{\partial t}
\]
We then take 
$\Pi$ times the $\Theta$ equation, $\partial p / \partial s$ times the
$\phi$ equations to derive
\[
 \dfrac{\partial }{\partial t} I  = 
 - c_p  \Pi \grad_s\cdot \left( \Theta \uvec \right)   
 - c_p  \Pi  \dfrac{\partial}  {\partial s}  \left( \Theta \dot s \right)
+ \frac{\partial p}{\partial s}   \uvec \cdot \grad_s \phi 
+ \frac{\partial p}{\partial s} \left(  \dot s   \frac{\partial \phi}{\partial s} \right)
- \frac{\partial p}{\partial s}  gw
\label{E:IEanalytic}
\]
The first, third and fifth terms represent transfer to kinetic energy
and will cancel after integration by parts
with similar terms in the kinetic energy
equation.  
In order for the terms involving $\dot s$  to cancel, we require
 \begin{equation}
  \int  c_p  \Pi  \dfrac{\partial}  {\partial s}  \left( \Theta \dot s \right) \, ds
 =
 \int  \frac{\partial p}{\partial s} \left(  \dot s   \frac{\partial \phi}{\partial s} \right) \, ds.
 \label{E:thetabarcondition}
 \end{equation}
this relation is needed for the advection terms to
be energetically neutral.  It can be achieved if
 care is taken in how the vertical advection of $\Theta$ is defined, or
 if we are using a vertically Lagrangian coordinate where these terms
 do not appear.

\subsection{Kinetic energy}
\label{S:KE}
For kinetic energy, take the sum of $\dpids \uvec$ dotted with the $\uvec$ equation,
$\dpids w$ times the $w$ equation 
and $(\uvec^2+ w^2)/2$ times the continuity equation.  First we consider only the
terms related to advection, which should be energetically neutral.  After removing two  terms that cancel pointwise, we are left with
\begin{multline}
 \dfrac{\partial }{\partial t} K = 
- \frac12 \dpids \uvec \cdot  \grad_s  \uvec^2 
- \frac12 \uvec^2  \nabla_s \cdot \dpids \uvec 
- \dpids  \uvec \cdot (\dot s  \frac{\partial \uvec}{\partial s})
- \frac12 \uvec^2        \frac{\partial}{\partial s}(\dpids \dot s)
\\
- \dpids  \uvec \cdot \grad_s \frac12 w^2  
- \frac12  w^2     \nabla_s \cdot \dpids \uvec 
- \dpids  w (\dot s \frac{\partial w}{\partial s})
- \frac12  w^2           \frac{\partial}{\partial s}(\dpids \dot s)
+ \dots
\\
\label{E:KEbudget1}
\end{multline}
The first, second, fifth and sixth terms in the above will cancel when integrating over the sphere via
integration by parts.  The third, forth, seventh and eight terms cancel in the continuum but rely on additional properties that are not obeyed by general mimetic discretizations. They are responsible for a small lack of conservation in \citeA[\S~4.4]{Dubos2015}. 
These terms will cancel via special integral properties satisfied by the SB81 discretization \cite{SB81,Taylor11}.  In \S~\ref{S:SB81identities} we show that
this result can be obtained more directly by noting that
SB81 also has a discrete derivative product rule.

Now consider the remaining terms involving the pressure gradient
\begin{equation}
 \dfrac{\partial }{\partial t} K = \dots 
 -gw \dpids
 + gw \dpids  \mu 
 - \uvec \dpids  \left( c_p  \theta_v \grad_s \Pi   + \mu \grad_s \phi \right).
\label{E:KEbudget2}
\end{equation}
The first term balances the identical term in the potential energy equation
and the remaining terms balance with similar terms in the
internal energy equation.

\subsection{Energy consistency}
\label{S:energyconsistency}
For energy consistency, we require that the discrete equations
reproduce three key aspects of the continuum equations: (1) the
transport terms are energetically neutral, meaning that they vanish to
machine precision when integrated over the domain, (2) the
diagnostic equations for the change in total potential, internal and
kinetic energies hold in exact time integration and (3) the transfer
terms in these diagnostic equations cancel to machine precision.
These three requirements ensure that the discretization is not introducing
any spurious sources of energy and that total energy will be conserved
in exact time integration.

We denote integration over the entire domain with a double integral, representing
integration in the vertical with respect to $ds$ and integration in the horizontal
with respect to an area metric $dA$.  The energy transfer terms in the budget
equations are denoted
\begin{equation}
  T_1   =   \iint c_p  \Theta \uvec \cdot  \grad_s \Pi     \qquad
\qquad
T_2 =  \iint gw \dpids
\qquad
T_3   = 
 \iint \left( \uvec \cdot \grad_s \phi -gw \right) \frac{\partial p}{\partial s}
\end{equation}
\begin{equation}
S_1    = 
  \iint c_p  \Pi \grad_s\cdot \left( \Theta \uvec \right)   
\qquad
\end{equation}
After integrating the energy budgets derived in the previous subsection, we
obtain 
\begin{align}
  \frac{\partial }{\partial t} \iint P &=  T_2
  \label{E:dPdt}
  \\
  \frac{\partial }{\partial t} \iint I &=  -S_1+T_3
  \label{E:dIdt}
  \\
  \frac{\partial }{\partial t} \iint K &=  -T_1 - T_2 -T_3
      \label{E:dKdt}
\end{align}
Energy conservation requires $T_1=-S_1$, which holds via integration-by-parts.

\section{Vertical discretization}
\label{S:vert_disc}
\begin{figure}
\begin{center}
\includegraphics[width=3.5in,trim=1.9in 3.75in 1.75in 3.75in,clip]{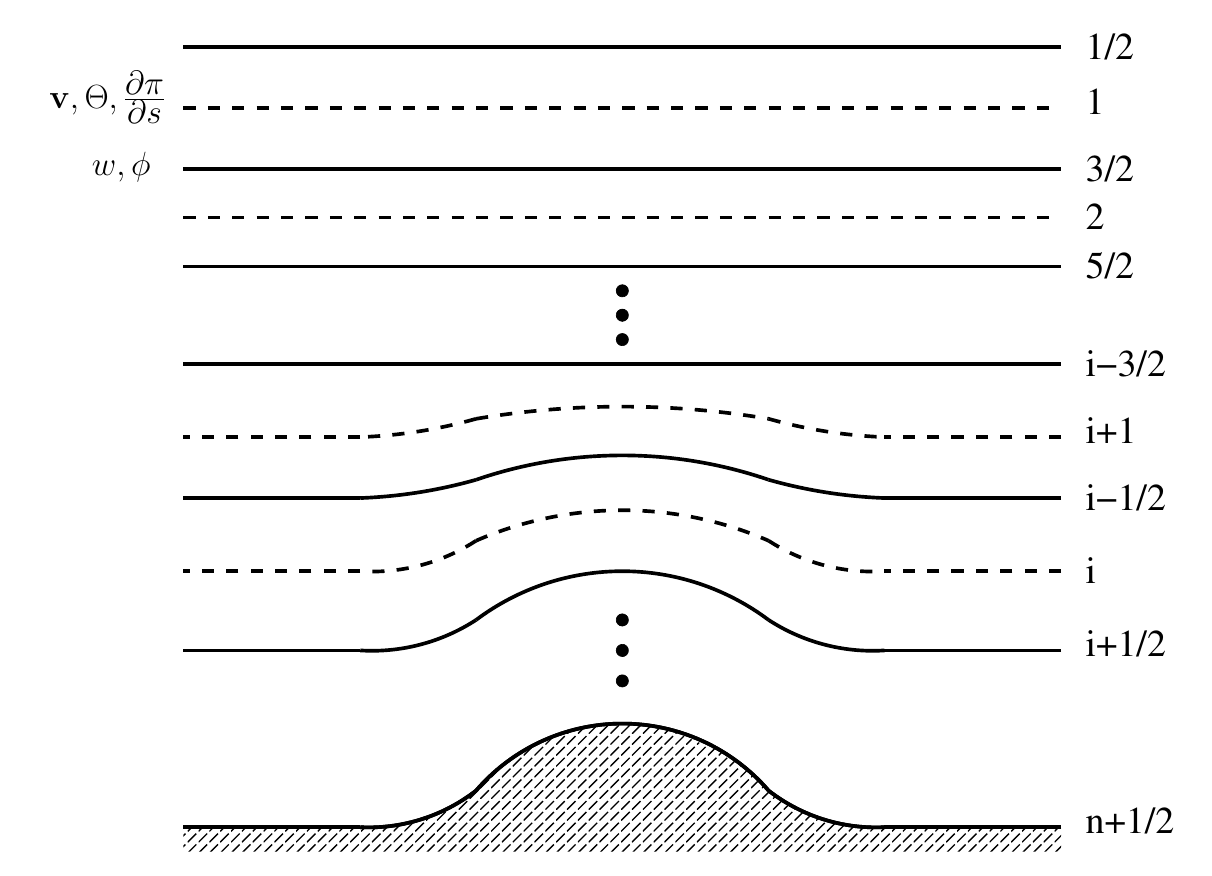}
\end{center}
\caption{
The terrain following $s$-coordinate levels
and level indexing.  There are
$n$ midpoint levels indexed by
$i=1,2,\dots,n$ and $n+1$ interface levels indexed by $i+\half,  i=0,1,\dots,n$.  
 }
\label{F:etacoords}
\end{figure}

We use a \citeA{Lorenz1960} staggering where the prognostic variables
are defined at either level interfaces or midpoints.  These levels, their indexing
and the staggering of the prognostic variables is shown in Fig.~\ref{F:etacoords}. The placement of the prognostic variables then leads to the following natural placement
for the remaining variables,
\begin{align}
\text{midpoint quantities:}&\qquad
\uvec, \Theta, \theta_v, \dpids,  p, \Pi, \frac{\partial \phi}{\partial s}
\\
\text{interface quantities:} &\qquad \phi, w, \pi, \dot s, \dot S , \mu, \frac{\partial p}{\partial s}  
\end{align}
where for convenience we have introduced $\dot S = \dot s \dpids$.
The formulas for $\dot S$ and $\mu$ mix midpoint and interface quantities,
and the averaging necessary for energy consistency will be given in the
next section. For the level positions, we require the midpoints to be centered within
the interfaces
\[
 s_i = \frac12( s_{i+\half} + s_{i-\half} )
\]
and denote the level thickness by $\DS$ with
\[
\DS_i = s_{i+\half} - s_{i-\half}
\qquad
\DS_{i+\half} = s_{i+1} - s_{i}
\]
and at the boundaries, $\DS_\text{n+\half}=\DS_n$ and  $\DS_{\half}=\DS_1$.
We note that the coordinate labeling by $s$ is arbitrary.
The solution will depend on the physical position of each level, but
not its assigned $s$-value.  
For simplicity we could relabel $s$ so that  $\Delta s = 1$ for all levels,
as is often done in numerical implementations.

Below we will define the averaging and differential operators using
a generic $\phi_{i+\half}$ variable 
given at interface points and a generic 
$p_i$ variable given at midpoints with boundary conditions
at the model top $p_{\half}$  and surface $p_{n+\half}$.
We use standard $\DS$ weighted averaging to map between
interface and midpoint quantities.  Because of the centered
midpoint levels, to average an interface quantities to midpoints
we use
\begin{equation}
 \bar \phi_i = \frac12 \left( \phi_{i+\half} + \phi_{i-\half} \right)
\label{E:avei2m}
\end{equation}
and to average (and extrapolate) midpoint quantities
to interfaces, we use
\begin{equation}
  \bar p_{i
    +\half} = \frac{(p \DS )_{i+1} + (p \DS )_i }{2 \DS_{i+\half}}
\qquad
\bar p_{\half} = p_1
\qquad
\bar p_{n+\half} = p_n
  \label{E:avem2i}
\end{equation}
This simple extrapolation formula is necessary to ensure a 
discrete averaging-by-parts identity holds. 
Vertical derivatives are defined by
\begin{equation}
\left(\frac{\partial p}{\partial s}\right)_{i+\half} = \frac{p_{i+1}-p_i}{\DS_{i+\half}}
\qquad
\left(\frac{\partial \phi}{\partial s}\right)_{i} = \frac{\phi_{i+\half}-\phi_{i-\half}}{\DS_{i}}
\label{E:vertderiv1}
\end{equation}
with the boundary conditions and one sided differencing used for 
computing derivatives at the surface and model top
\begin{equation}
\left(\frac{\partial p}{\partial s}\right)_{\half} =  \frac{  p_1 -  p_{\half}  }{\frac12 \DS_{i+\half}}
\qquad
\left(\frac{\partial p}{\partial s}\right)_{n+\half} = 
\frac{  p_{n+\half} - p_{n} }{\frac12 \DS_{n}}
\label{E:vertderiv2}
\end{equation}
Note that for the $\partial/\partial s$ operator on midpoint
quantities, these boundary conditions must be specified as part
of the discretization.  

Vertical integrals are approximated via quadrature using a midpoint rule.
For midpoint quantities,
\begin{equation}
  \int p \, ds \approx  \sum_{i=1}^n   p_i  \DS_i
  \label{E:quad-midpoint}
\end{equation}
and at interfaces
\begin{equation}
  \int \phi \, ds \approx   {\sum_{i=0}^n}'  \phi_{i+\half} \DS_{i+\half}
  \label{E:quad-interface}  
\end{equation}
where the prime on the sum $\sum'$ is used to denote that
we multiply the integrand by $\frac12$ at the end points $i=0$ and $i=n$.  

\subsection{Discrete operator identities}
\label{S:SB81identities}

With these definitions, our discrete operators obey integration
and averaging by parts identities as well as a discrete
derivative product rule.  For interface quantities, averaging
commutes with differentiation
\begin{equation}
\overline{ \left( 
 \frac{\partial \phi}{\partial s} 
 \right)}
=
\frac{\partial \bar\phi}{\partial s}
\label{E:bards}
\end{equation}
with the boundary values of $\phi$ used for the
$\partial \bar \phi / \partial s$ boundary conditions.
We have a discrete averaging-by-parts identity
 \begin{equation}
{\sum_{i=0}^n}'  \left(\bar p\right)_{i+\half} \, \phi_{i+\half}  \, \DS_{i+\half}
= 
\sum_{i=1}^n  p_i
\left(\bar\phi\right)_i   \, \DS_i
\label{E:avebyparts}
\end{equation}
and a discrete integration-by-parts identity 
\[
{\sum_{i=0}^n}'  \phi_{i+\half}  \frac{\partial p}{\partial s} {\DS_{i+\half}} 
+
\sum_{i=1}^n  \frac{\partial \phi}{\partial s} p_i  \DS_i
=   \phi_{n+\half} p_{n+\half}
- \phi_\half p_\half
\]
Most importantly for our formulation, and difficult to obtain for general
discretizations \cite{Salmon07}, is the pointwise product
rule.  For interface quantities $a$ and $b$, we have
\begin{equation}
\left( \frac{\partial   }{\partial s} (a b) \right)_{i}
=
\bar b \left(\frac{\partial a }{\partial s}\right)
+
\bar a \left(\frac{\partial b }{\partial s}\right)
\label{E:productrule1}
\end{equation}
For midpoint quantities $c$ and $d$ we have a similar formula,
but we have to use unweighted averaging,
\begin{equation}
\left( \frac{\partial   }{\partial s} (c d)\right)_{i+\half}
=
\overline{ \left(\frac{d}{\DS_i} \right)} \DS_{i+\half}    \left(\frac{\partial c }{\partial s}\right)
+
\overline{ \left(\frac{c}{\DS_i} \right)} \DS_{i+\half}   \left(\frac{\partial d }{\partial s}\right)
\label{E:productrule2}
\end{equation}
where we have extra factors of $\DS$ so that the unweighted average
of midpoint quantities 
is written in terms of the weighted average defined in \eqref{E:avem2i}.

\subsection{SB81 vertical advection}
\label{S:SB81advection}
The SB81 vertical advection for operator $\dot s \partial / \partial s$ 
is carefully constructed for midpoint quantities in order to
conserve energy.
We denote this operator with square brackets $[\cdot]$ and write
it in terms of interface quantitiy $\dot S$.
For midpoint quantities,
\begin{equation}
\dpids \left[ \dot s \frac{\partial p}{\partial s} \right]_{i} 
=
\frac{1}{\DS_i} \overline{ \left( 
   \dot S \frac{\partial p}{\partial s} \DS_{i+\half} 
\right)}
=
\left( 
 \frac{\partial }{\partial s}  \left( \dot S 
   \overline{ \frac{ p}{ \DS_i}} \DS_{i+\half}  \right)
- p \frac{\partial  }{\partial s} \dot S.  
\right)
\label{E:SB81vertadv1}
\end{equation}
For interface quantities, we introduce a generalization of this
operator defined by
\begin{equation}
\overline\dpids \left[ \dot s \frac{\partial \phi}{\partial s} \right]_{i+\half} 
=
 \overline{ \left( 
\overline{\dot S} \frac{\partial \phi}{\partial s} 
\right)}
=
\left( 
 \frac{\partial }{\partial s}  \left( \overline{\dot S} \bar \phi \right)
- \phi \frac{\partial }{\partial s}  \overline{\dot S}
\right)
\label{E:SB81vertadv2}
\end{equation}
where we have averaged the midpoint quantity $\dpids$ to interfaces.
This generalization appears to contain double averaging but it can also be written in an equivalent form without averages of averages.

\section{Lorenz staggered formulation}

We now write the equations in the vertically staggered
discretization using the bar notation
for averaging between midpoints and interfaces using 
\eqref{E:avem2i} or \eqref{E:avei2m}, $\partial/\partial s$
representing \eqref{E:vertderiv1} and  \eqref{E:vertderiv2}
and $[\cdot]$ representing the SB81 transport operator.
The discrete equations in their final energy consistent form
are 
\begin{align}
  \frac{\partial \uvec}{\partial t} + (\nabla_s \times  \uvec + f\khat)  \times \uvec 
  +\frac12 \grad_s \left( \uvec^2 + \overline{ w^2} \right)
   - \overline{w \grad_s w} 
   + \left[\dot s \frac{\partial }{\partial s} \right] \uvec
\label{E:uequation-discrete}
\\
+
c_p  \theta_v  \grad_s \Pi 
+\overline{ \mu \grad_s \phi  }
 &= 0   
\\
  \frac{\partial w}{\partial t} +
  \tilde \uvec \cdot \grad_s w +
  \left[\dot s \frac{\partial }{\partial s} \right] w
+g
-g \mu 
  &= 0
\\
\frac{\partial \phi}{\partial t} +
\tilde \uvec 
\cdot \grad_s \phi +
  \dot s \frac{\partial\bar\phi }{\partial s} 
 - gw &= 0
\\
\dfrac{\partial  }{\partial t}   \Theta 
 + \grad_s\cdot \left(  \Theta \uvec \right)   + 
  \dfrac{\partial}  {\partial s}  \left( \tilde\theta_v \dot S \right)   &= 0
  \label{E:thetaequation-discrete}
\\
\dfrac{\partial }{\partial t} \left( \dpids \right)  + \nabla_s \cdot  \left( \dpids  \uvec \right)    + 
\dfrac{\partial }{\partial s}  \dot S  &= 0
\label{E:rhoequation-discrete}
\end{align}
with boundary condition
\[
\dot s_{\half} = \dot s_{n+\half} = 0.
\]
With this staggering, $p$ and $\Pi$ at their midpoint locations
are determined from \eqref{E:EOS} without averaging,
\[
\frac{\partial \phi}{\partial s}  = - R \Theta \frac{\Pi}{p}
\qquad
\Theta = \dpids  \theta_v 
\]
and averaging is needed for vertical momentum and $\mu$,
\[
\dot S = \overline\dpids \dot s
\qquad
\mu_{i+\half} = \left( \frac{\partial p}{\partial s}  \Big/  \overline \dpids \right).
\]
In order for the transport terms to be energetically neutral,
we need to use special averaging denoted by a tilde 
for two variables:
\begin{equation}
\tilde\uvec_{i+\half} =   \overline{\left( \dpids \uvec \right)} \Big/ \overline\dpids 
\qquad
    \left(\tilde\theta_v\right)_{i+\half} =
- \frac{\mu}{c_p}   \overline{ \frac{\partial \phi}{\partial s} } \Big/  \dfrac{\partial \Pi  }  {\partial s} 
\label{E:tildedefinition}
\end{equation}
For $\tilde \uvec$ this is the density weighted average.
For $\bar \theta_v$ this represents a version of \eqref{E:EOS} at interfaces.
This term
will not be present in the vertically Lagrangian coordinate.
In the Eulerian coordinate, \eqref{E:EOS}
at interfaces is thus treated quite differently then \eqref{E:EOS} at midpoints.
As an alternative one can eliminate \eqref{E:EOS} altogether
by introducing an additional prognostic variable for $p$ or $\Pi$
\cite{Gassmann2013}.

For the mass coordinate, we use a diagnostic integral for $\dot S$, $i=1 \dots n-1$
\[
\dot S_{i+\half} = 
 B_{n+\half} \sum_{i=1}^n  \grad_s \cdot \left( \dpids \uvec \right) \,\DS_i
- \sum_{i=1}^i  \grad_s \cdot  \left( \dpids \uvec \right)  \,\DS_i
\]
and then compute $\dot s$.  For the Eulerian height coordinate, we replace the prognostic
equation for $\phi$ with the diagnostic equation
\[
  \dot s \frac{\partial\bar\phi }{\partial s} 
 = -\tilde \uvec \cdot \grad_s \phi  + gw
\]
and then compute $\dot S$.  
In the vertically Lagrangian case, $\dot S = \dot s = 0$ and all terms
involving these variables are dropped.  In the vertically Lagrangian
height coordinate, we retain the $\phi$ equation with $\phi$ at the
surface and model top held fixed.  

We require that the equations for $w$ and $\phi$ hold at all
interfaces including the model top and surface.  This requirement is
used to deduce the consistent boundary conditions for $w$ and $p$.  At
the surface, $\phi$ is given for both mass and height coordinates and
the $\phi$ equation is satisfied by taking $g w = \tilde \uvec \cdot
\grad_s \phi_s $.  This allows us to determine $\partial w / \partial
t$ at the surface in terms of $\partial \uvec / \partial t$,
which in turn determines the pressure (and equivalently $\mu$) at the
surface using the prognostic equation for $w$.  We omit this algebra
here as it depends on details of the spatial discretization.
For our purposes, we need only know that $p$ and $\mu$ are defined
at the surface so that the prognostic $w$ equation holds at the surface.

At the model top, for the height coordinate the procedure is identical and
$w$ and $\mu$ are defined so that the $w$ and $\phi$ equations hold at the
top most interface.  For the mass coordinate, we have a pressure boundary condition
and a free surface at the model top.  In this case, $\mu$ is well defined and $w$ and
$\phi$ are evolved at the model top with their prognostic equations.

We note that we used the SB81 operator for vertical advection of
$w$ and the more traditional discretization for vertical advection of $\phi$.
The SB81 operator is necessary for energy consistency in the momentum equations,
but it does not vanish at the surface.  The traditional discretization does
vanish at the surface since $\dot s = 0$, and its use in the $\phi$ equation
results in the more natural $g w = \tilde \uvec \cdot \grad_s \phi_s $
boundary condition for $w$.

\subsection{Discrete energy conservation}

We now repeat the energy budget calculations done for the
continuum formulation in \S~\ref{S:EnergyConservation}.
For our horizontal discretization we will be using the mimetic
collocated spectral element method from \citeA{TayFour10}.
In this case, the horizontal terms which cancel after integration
by parts will cancel in the discrete case as well and those
calculations are not repeated here.  For a horizontally staggered
case, such as \citeA{thuburn2009numerical}, more care must be taken
to ensure the horizontal integration by parts properties still hold
when mixed with vertical staggering \cite{Gassmann2013}.

\subsection{Discrete potential energy}
We define discrete potential energy at midpoints
\[
{\sum_{i=1}^n} P  \DS_{i}  
  = \sum_{i=1}^n  \dpids  \bar\phi \DS_{i}
  = {\sum_{i=0}^n}' \overline\dpids  \phi \DS_{i+\half}
\]
For the budget, we derive the discrete analog of equations
\eqref{E:dPdt}.  We take the 
continuity equation, average to interfaces
and multiply by $\phi$, and then 
add in the $\phi$ equation multiplied by
$\overline\dpids$. Assuming exact time integration, we have 
\begin{multline}
{\sum_{i=1}^n} \frac{\partial P}{\partial t}  \DS_{i} 
=
- {\sum_{i=0}^n}'  \phi \nabla_s \cdot  \overline{\left( \dpids  \uvec \right)     }  \,\DS_{i+\half}
-  {\sum_{i=0}^n}'  \phi \overline{\dfrac{\partial  \dot S}{\partial s}  } \,\DS_{i+\half}
\\
- {\sum_{i=0}^n}' \overline\dpids  \tilde\uvec \cdot \grad_s \phi  \,\DS_{i+\half}
-  {\sum_{i=0}^n}' \overline\dpids  \dot s \overline{\frac{\partial \phi}{\partial s}  }  \,\DS_{i+\half}
+ {\sum_{i=0}^n}' \overline\dpids gw  \,\DS_{i+\half}.
\end{multline}
Similar to the continuum equations, four of these terms will cancel when integrated
over the domain after integration by parts and making use of
the definition of $\tilde \uvec$. The last term, after averaging by parts, becomes,
\begin{equation}
  {\sum_{i=1}^n} \frac{\partial P}{\partial t}  \DS_{i} 
= \dots + 
 {\sum_{i=1}^n} \dpids g \bar w  \,\DS_{i}.
\end{equation}

\subsection{Discrete internal energy}
We use a discrete internal energy 
\[
\sum_{i=1}^n I \DS_i = 
\sum_{i=1}^n  c_p \Theta \Pi \DS_i 
+{\sum_{i=1}^n} \left( \frac{\partial\phi}{\partial s} \right)_i  p_i  \DS_i
+ \hat p \phi_{\half} 
\]
with $\hat p = p_\half$ in the mass coordinate and $\hat p = 0$ in the height
coordinate.  Differentiating with respect to time, 
and taking $\partial \Pi / \partial t = (R/c_p) (\Pi / p) \partial p/ \partial t$
which will hold in exact time integration as well as in the
the limit of small $\Delta t$, and using that $(\partial \hat p/ \partial t) = 0$
in both coordinate systems, 
\[
\sum_{i=1}^n \frac{\partial I}{\partial t} \DS_i = 
\sum_{i=1}^n  c_p  \frac{\partial\Theta}{\partial t} \Pi \DS_i 
+{\sum_{i=1}^n}
\frac{\partial}{\partial s} \left( \frac{\partial \phi}{\partial t} \right)
p_i  \DS_i 
+ \hat p \frac{\partial \phi_\half}{\partial t} 
\]
where we have removed the terms
\[
\sum_{i=1}^n  R \Theta \frac{\Pi}{p} \frac{\partial p}{\partial t} \DS_i 
+{\sum_{i=1}^n} \frac{\partial \phi}{\partial s} \frac{\partial p}{\partial t}
\DS_i 
=0
\]
since our discrete EOS is 
$
\frac{\partial \phi}{\partial s}  = -R ( \Theta \Pi / p)_i.
$
We then apply integration by parts to obtain
\[
\sum_{i=1}^n \frac{\partial I}{\partial t} \DS_i = 
\sum_{i=1}^n  c_p  \frac{\partial \Theta}{\partial t} \Pi \DS_i 
-{\sum_{i=0}^n}' \frac{\partial \phi}{\partial t} \frac{\partial p}{\partial s}
\DS_{i+\half}. 
\]
Inserting the discrete prognostic equations,
\begin{multline}
  \sum_{i=1}^n \frac{\partial I}{\partial t} \,\DS_i = 
 -  \sum_{i=1}^n c_p  \Pi_i \grad_s\cdot \left( \Theta \uvec \right)  \,\DS_i
 - \sum_{i=1}^n c_p  \Pi_i   \dfrac{\partial}  {\partial s}  \left(  \tilde\theta_v  \dot S\right)    \,\DS_i
\\
+ {\sum_{i=0}^n}' \frac{\partial p}{\partial s}  \tilde \uvec \cdot \grad_s \phi  \,\DS_{i+\half}
+ {\sum_{i=0}^n}' \frac{\partial p}{\partial s}   \dot s \overline{  \frac{\partial \phi}{\partial s} } \,\DS_{i+\half}
\\
- {\sum_{i=0}^n}' \frac{\partial p}{\partial s}  gw   \,\DS_{i+\half}
\label{E:IElorenz}
\end{multline}
The first, third and fifth terms make up the transfer terms which will balance
similar terms in the kinetic energy equation.  
For the vertical advection terms, we integrate one term by parts and
use the definition of $\dot S$ to obtain
\[
\sum_{i=1}^n \frac{\partial I}{\partial t} \,\DS_i =
\dots +
  {\sum_{i=0}^n}' c_p  \dfrac{\partial \Pi  }  {\partial s} 
  (\tilde\theta_v)_{i+\half}  \dot S 
  \DS_{i+\half}
  + {\sum_{i=0}^n}' \mu  \dot S      \overline{ \frac{\partial \phi}{\partial s} }
    \DS_{i+\half}.
\]
Ensuring that these two terms cancel is the reason behind our
complex formula for $\tilde \theta_v$.

\subsection{Discrete kinetic energy}
Starting with $K$ at midpoints,
\[
K_i = \frac12 \dpids \uvec^2 + \frac12 \dpids \overline{w^2}
\]
summing and applying averaging by parts we have
\begin{equation}
\sum_{i=1}^n K_i \,\DS_i = 
\frac12 \sum_{i=1}^n   \dpids \uvec^2  \,\DS_i
+\frac12  {\sum_{i=0}^n}'  \overline\dpids {w^2}  \,\DS_{i+\half}
\label{E:KELorenz}
\end{equation}
We then derive the discrete budget following \S~\ref{S:KE}.
We dot the discrete equation for $\uvec$ by $\dpids \uvec$ and sum with
$\dpids w$ times the $w$ equation
and $(\uvec^2+ w^2)/2$ times the continuity equation.
Starting with just the advection terms from \eqref{E:KEbudget1},
we drop four of the terms which cancel after integration by parts
in the horizontal and two terms cancel after integration by
parts in the vertical, using the SB81 operator for $\dot s \partial \uvec / \partial s$
and the product rule identity \eqref{E:productrule2}.
The remaining advection terms are 
\begin{multline}
\sum_{i=1}^n \frac{\partial K_i}{\partial t} \,\DS_i =
- {\sum_{i=0}^n}'  \overline\dpids w \tilde\uvec \cdot \grad_s w   \,\DS_{i+\half} 
+ {\sum_{i=1}^n} \dpids \uvec \cdot \overline{w \grad_s w}
\\
- {\sum_{i=0}^n}' \overline\dpids  w \left[\dot s \frac{\partial w}{\partial s}\right]_{i+\half}  \,\DS_{i+\half} 
- {\sum_{i=1}^n} \overline{\frac12  w^2}           \frac{\partial}{\partial s} \dot S  \,\DS_{i}
+\dots
\end{multline}
Using the definition of $\tilde \uvec$ from \eqref{E:tildedefinition} the first two terms
cancel.  Expanding the  SB81 operator and applying integration-by-parts we are
left with
\[
\sum_{i=1}^n \frac{\partial K_i}{\partial t} \,\DS_i =
- {\sum_{i=1}^n}   \overline w   \frac{\partial w}{\partial s} 
\overline{\dot S} 
   \,\DS_{i} 
- {\sum_{i=1}^n}  \overline{\left( {\frac12  w^2}     \right)}
    \frac{\partial}{\partial s}   \dot S \,\DS_{i}  +\dots
\]
To show these two terms cancel, we apply the SB81 product rule,
average by parts, apply \eqref{E:bards} and then integrate
by parts making use of $\dot S =0$ at the boundaries,
\begin{multline}
- {\sum_{i=1}^n}   \overline w   \frac{\partial w}{\partial s} 
\overline{\dot S} 
   \,\DS_{i} 
   =
   - {\sum_{i=1}^n}     \frac{\partial }{\partial s} \left(\frac12 w^2 \right)
\overline{\dot S} 
   \,\DS_{i} 
   =
   - {\sum_{i=0}^n}'  \overline{ \frac{\partial }{\partial s} \left(\frac12 w^2 \right)}
{\dot S} 
\,\DS_{i+\half}
\\
= - {\sum_{i=0}^n}' \frac{\partial}{\partial s} \overline{\left( {\frac12  w^2}     \right)}
 {      \dot S }\,\DS_{i+\half}  
=  {\sum_{i=1}^n} \overline{\left( {\frac12  w^2}     \right)}
\frac{\partial}{\partial s}  {      \dot S }\,\DS_{i}
\end{multline}

Finally, we are left with the discrete version of the transfer terms
from \eqref{E:KEbudget2},
\begin{multline}
{\sum_{i=1}^n} \frac{\partial K_i}{\partial t} \,\DS_i = \dots
 - {\sum_{i=0}^n}' gw \overline\dpids  \,\DS_{i+\half} 
 + {\sum_{i=0}^n}' gw \overline\dpids  \mu
  \,\DS_{i+\half} 
\\
- {\sum_{i=1}^n} \uvec \dpids \cdot \left(
c_p  \theta_v  \grad_s \Pi 
+\overline{ \mu \grad_s \phi  }
 \right)
 \,\DS_{i} 
\end{multline}

\subsection{Discrete energy consistency}
We now write the discrete integral energy budgets.
Let the sum over $i$ represent the vertical integral as above.
In the horizontal, we let the sum over  $j$ represent the
discrete quadrature for integration over spherical $s$ surfaces
associated with the mimetic horizontal discretization.

\begin{equation}
  T_1   =   \sum_j {\sum_{i=1}^n}  c_p  \Theta \uvec \cdot  \grad_s \Pi \DS_{i}
  \qquad
T_2 =  \sum_j {\sum_{i=0}^n}' gw \overline\dpids \DS_{i+\half}
\qquad
\end{equation}
\begin{equation}
T_3   = 
\sum_j  {\sum_{i=1}^n} \uvec \dpids \cdot \overline{\mu \grad_s \phi}  \DS_{i}
- \sum_j {\sum_{i=0}^n}' gw  \frac{\partial p}{\partial s} \DS_{i+\half}
\end{equation}
\begin{equation}
S_1    = 
  \sum_j {\sum_{i=1}^n} c_p  \Pi \grad_s\cdot \left( \Theta \uvec \right)   \DS_{i}
  \qquad
  S_2 =
\sum_j {\sum_{i=1}^n} g \bar w \dpids \DS_{i}
\end{equation}
\begin{equation}
S_3   = 
\sum_j  {\sum_{i=1}^n} \frac{\partial p}{\partial s} \tilde \uvec \cdot \grad_s \phi   \DS_{i}
- \sum_j {\sum_{i=0}^n}' gw  \frac{\partial p}{\partial s} \DS_{i+\half}
\end{equation}

After integrating the energy budgets derived in the previous subsection, we
obtain 
\begin{align}
\frac{\partial }{\partial t}  \sum_j \sum_{i=1}^n P  &=  S_2
  \label{E:dPdt-discrete}
  \\
\frac{\partial }{\partial t} \sum_j \sum_{i=1}^n I &=  -S_1 + S_3
  \label{E:dIdt-discrete}
  \\
\frac{\partial }{\partial t} \sum_j \sum_{i=1}^n K &=  -T_1 - T_2 -T_3.
      \label{E:dKdt-discrete}
\end{align}
To obtain energy conservation, we note that $S_1=-T_1$ after integration by parts and
$S_2=T_2$ after averaging-by-parts.  To show $S_3=T_3$ we average by parts and then expand
$\mu$ and $\tilde \uvec$.  Ensuring $S_3=T_3$ is the motivation behind the averaging of
$\mu \grad_s \phi$ in the momentum equation.  This is the only term in our formulation
that contains double averaging (since $\mu$ contains an average of $\dpids$).
This is unavoidable given our definition of $\tilde \uvec$ needed to obtain
energy neutral transport.

\section{Numerical implementation}

We have implemented the discretization described here in an open
source version of HOMME.  
HOMME previously contained a shallow water and hydrostatic 
primitive equation model and we use HOMME-NH to refer to the new
nonhydrostatic model.  For the horizontal directions,
HOMME-NH relies on the existing HOMME infrastructure including
the mimetic spectral element discretization, monotone and
conservative tracer transport \cite{Guba2014}, and support for
fully unstructured arbitrary quadrilateral grids \cite{Guba2014tensor}.  
Here we will be running the
model on the quasi-uniform equi-angle cubed
sphere grid \cite{Ronchi96}.  Within each spectral element we use 4th
order accurate polynomial basis functions of degree $p=3$.  

HOMME-NH implements \eqref{E:uequation-discrete}--\eqref{E:rhoequation-discrete}
using method of lines.
The horizontal divergence, gradient
and curl operators on $s$-surfaces are evaluated with existing HOMME
subroutines, and the vertical terms are computed as written.
All the tendency terms
are computed locally within each element and then we apply a
standard finite element assembly step.  Because of the spectral
element method's diagonal mass matrix, this assembly
is explicit and can be computed in each element with only
nearest neighbor information.   For our initial HOMME-NH implementation we
use the mass coordinate for with both Eulerian and vertically Lagrangian
options.  For the vertically Lagrangian option, we typically remap back to
the initial reference levels every three timesteps using 
parabolic splines \cite{zerroukat06}.  

For timestepping, we use a HEVI approach \cite{Satoh2002} discretized with a Runge-Kutta IMEX method.
The HEVI/IMEX combination has proven attractive for atmospheric models, e.g. 
\citeA{UllrichJablonowski2012,GKC2013,Weller2013,Lock2014}.
HOMME-NH supports a variety of IMEX methods from
the ARKODE package \cite{Gardner2017,ARKODE2018} of the SUNDIALS library
\cite{hindmarsh2005sundials}.
Here we use a 5 stage 2nd order accurate high-CFL method derived in 
\citeA{Steyer2019}. 
In the HEVI approach, the vertically implicit component can be solved independently in each vertical column.  We use Newton iterations with the Jacobian computed analytically and inverted by
direct LU factorization.  We use a time-split approach to couple
\eqref{E:uequation-discrete}--\eqref{E:rhoequation-discrete} with
hyperviscosity and physical parameterizations following 
\citeA{Dennis12}.  

HOMME-NH is further described in \cite{ullrich2017dcmip} and results from
a splitting supercell thunderstorm test case are given in \cite{Zarzycki2019}.

\section{Numerical results}

We use HOMME-NH to verify the three energy consistent properties from 
\S~\ref{S:energyconsistency}.  Properties (1) and (3) depend on the
mimetic properties of the discretization.  These are verified
through unit tests for all the terms in the energy budgets
which should cancel when integrated over the domain. These tests are
useful for ensuring all differential operators and quadrature rules
are properly coded and are not presented here.   To establish
property (2) we verify that equations
\eqref{E:dPdt-discrete}-\eqref{E:dKdt-discrete} are satisfied to
within time truncation error when running HOMME-NH without forcing or
explicit dissipation.  
Since these equations would be exact in exact time integration,
their residual is a measure of the dissipation introduced by the time stepping algorithm.  
In an energy consistent discretization, this dissipation should converge to zero
with $\Delta t$.

We make use of a moist baroclinic instability test case from
the Dynamical Core Model Intercomparison Project (DCMIP)
\cite{ullrich2017dcmip,DCMIP16}.  This test case is based on
\citeA{ullrich2014} with the addition a prognostic variable for
specific humidity and the \citeA{Reed2012} simple physics package.
The physics includes large-scale condensation, surface fluxes and
boundary layer turbulence.  We use a low resolution grid with
$4^\circ$ average grid spacing at the Equator and 30 vertical levels.
This low resolution ensures energy conservation holds
even in the presense of large spatial truncation errors.
\begin{figure}
\begin{center}
  \includegraphics[width=2.5in,trim=.4in 3.0in 0.5in 2.9in,clip]{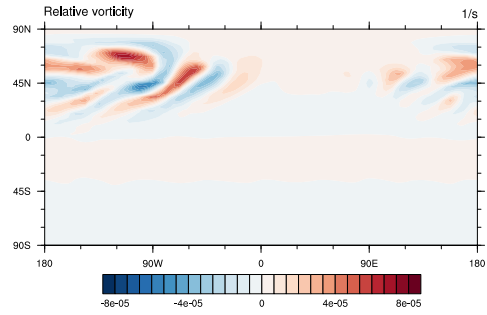}
  \end{center}
\caption{Relative vorticity on model level 21 (750 hPa)
  at day 15 in the DCMIP moist baroclincic test case.}
\label{F:vorticity}
\end{figure}
\begin{figure}
\begin{center}
  \includegraphics[width=2.5in,trim=.4in 3.2in 0.5in 2.9in,clip]{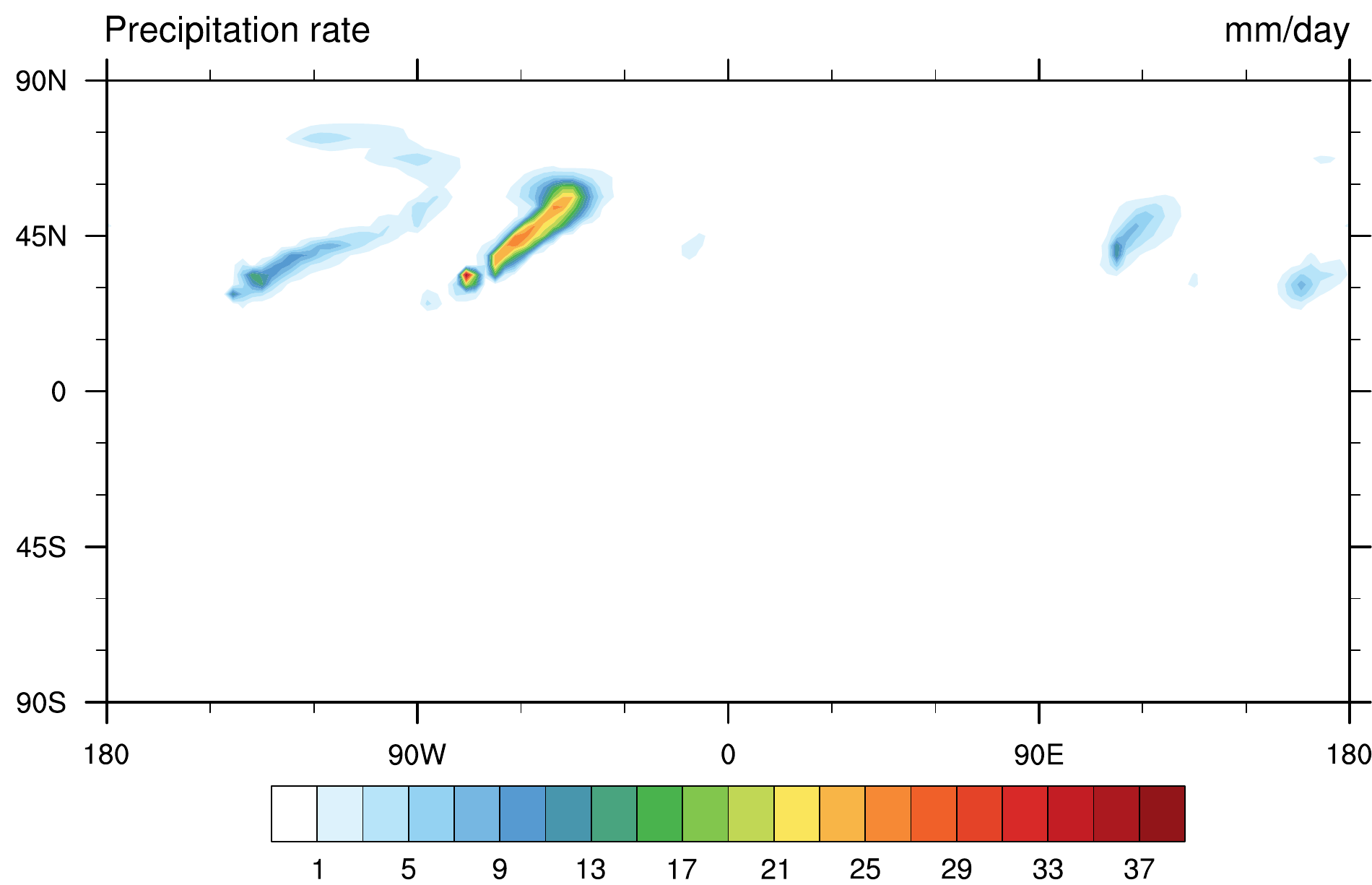}
  \end{center}
\caption{Precipiation rate at day 15 in the DCMIP moist baroclincic test case 
  generated by \citeA{Reed2012} simple physics.  }
\label{F:prec}
\end{figure}
To analyze the HOMME-NH discretization with a realistic flow containing energy in all resolved scales we first run the full model for 15 days.
The solution at $t_0=$15 days is shown in Figures~\ref{F:vorticity} and \ref{F:prec}.
We then disable the forcing and all explicit model
dissipation and run to time $t_1 = t_0 + $ 2 hours.
We use the two hour period $[t_0,t_1]$ to
study the HOMME-NH solution to the adiabatic component of the model, equations 
\eqref{E:uequation-discrete}--\eqref{E:rhoequation-discrete}.   
\begin{figure}
\begin{center}
  \includegraphics[width=3.5in,trim=0in 2.5in 0.5in 2.8in,clip]{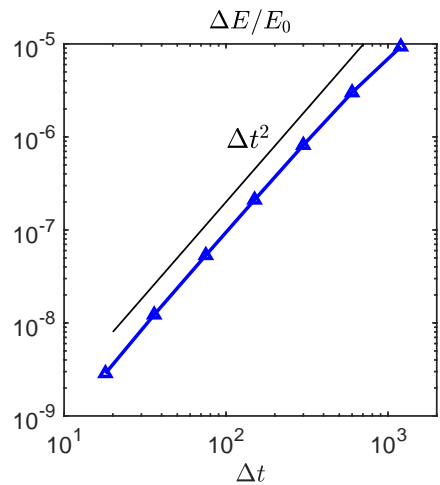}
  \end{center}
\caption{Relative change in total energy over two hours of model time.
  Energy is conserved to $O(\Delta t^2)$ with a second order Runge-Kutta
  method.  }
\label{F:tote}
\end{figure}
The energy conservation of the discretization as a function of $\Delta t$
is shown in Fig.\ref{F:tote}.  We plot the relative change in total energy
over $[t_0,t_1]$
\[
 \frac{E(t_1) - E(t_0)}{E(t_0)}
\]
for simulations with $\Delta t$ ranging from close to the CFL stability limit
(1200s) down to 18s.   With a second order accurate time stepping algorithm we expect
second order convergence in the energy, as shown.

\begin{figure}
\begin{center}
  \includegraphics[width=3.5in,trim=0in 2.7in 0.0in 2.7in,clip]{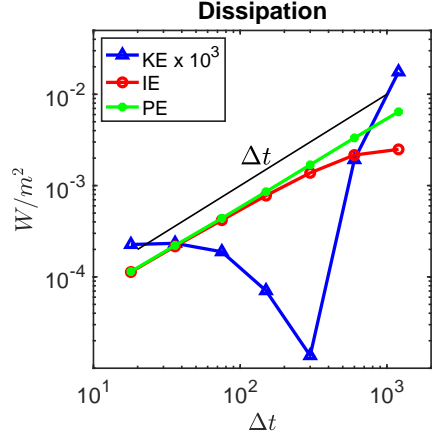}
  \end{center}
\caption{Residual error in the energy budget equations
  \eqref{E:dPdt-discrete}-\eqref{E:dKdt-discrete} when running
  without any additional dissipation.  This residual measures
  the dissipation from time truncation errors. The 
  potential and internal energy dissipation vanishes as $O(\Delta t$).
  The kinetic energy dissipation is three orders of magnitude smaller
  and vanishes at a faster rate until reaching machine precision.
}
\label{F:dissipation}
\end{figure}
Next we examine the residual
error for \eqref{E:dPdt-discrete}-\eqref{E:dKdt-discrete} instantaneously at time $t_1$ given
by
\begin{align}
  R_P &=  
  \sum_j \sum_{i=1}^n \frac{P(t_1+\Delta t)-P(t_1)}{\Delta t} -  S_2(t_1)
  \\
  R_I &=
  \sum_j \sum_{i=1}^n \frac{I(t_1+\Delta t)-I(t_1)}{\Delta t}  + S_1(t_1) - S_3(t_1)
  \\
  R_K &=
  \sum_j \sum_{i=1}^n \frac{K(t_1+\Delta t)-K(t_1)}{\Delta t}  +T_1(t_1) + T_2(t_1) +T_3(t_1).
\end{align}
The residuals $R_p, R_I$ and $R_k$ represent unphysical changes in
potential, internal and kinetic energy in HOMME-NH as a function of $\Delta t$.  
For simplicity we estimate $\partial/\partial t$ with a first order accurate
approximation so that
$R_p, R_I$ and $R_k$ can be computed in the code at each timestep with only
data from that timestep.  In Fig.\ref{F:dissipation} we plot $|R_P|,|R_I|$ and $|R_K|$,
with $\Delta t$ again ranging from
the CFL limit (1200s) to 18s.
From the figure we see that $R_p$ and $R_I$
decrease to zero as $O(\Delta t)$, as expected due to our first order
approximation.  The value of $R_k$ 
is 3 orders of magnitude smaller than $R_p$ and $R_I$, and decreases to
zero at an even faster rate until reaching machine precision.
For this problem, $R_K$ and $R_I$ are always negative, $R_P$ is always positive, 
$K=8.8\times10^5 \, J/m^2$, $I=1.8\times10^9 \, J/m^2$ and $P=7.2\times10^8 \, J/m^2$. 
At $\Delta t=300s$, $R_k = -1.4\times10^{-8} \, W/m^2$ corresponds 
to change in the 12'th digit of $K$.   
These results verify the HOMME-NH implementation of
\eqref{E:uequation-discrete}--\eqref{E:rhoequation-discrete},
confirming that the spatial discretization does not have
any spurious sources of energy and that the only changes in energy
are due to the timestepping algorithm.

Most simulations with HOMME-NH are run with timesteps close
to the CFL limit.  In this regime, we see that the dissipation from
time truncation is $\sim 0.02 \, \mbox{W/m}^2$ at day 15.
To compare this number with dissipation from the non-adiabatic
terms, we note that at day 15 before the explicit dissipation was turned off,
the kinetic energy dissipation from hyperviscsoity and the vertical remap
is $0.24 \, \mbox{W/m}^2$ and $0.06 \, \mbox{W/m}^2$ respectively.

\section{Conclusions}

We have presented a discretization of the nonhydrostatic equations in terrain
following coordinates with a Lorenz staggered vertical discretization.
The discretization is energy consistent in exact time integration when
coupled with a mimetic horizontal discretization, meaning the discrete
diagnostic equations for potential, internal and kinetic energy are
consistent and total energy is conserved to time truncation error. It
supports both mass and height coordinates, which differ only in the
top of model boundary condition and the diagnostic equation for $\dot
s$.  Vertical transport is kept in advective form, leading to unified
support for both Eulerian or vertically Lagrangian coordinates.
We use virtual potential temperature and prognose total mass to
support general moist equations of state.  The mass coordinate
formulation has been implemented in HOMME-NH, with energy
consistency verified using a moist baroclinic instability test case
and a Runge-Kutta IMEX timestepping algorithm.
In the design of HOMME-NH, we have found energy consistency to be a
useful guiding principal that enhances the stability of the discretization.
The ultimate goal is not exact energy conservation, but to preserve the Hamiltonian structure of the equations and ensure
all changes in energy are through explicitly added dissipative mechanisms.
In HOMME, these mechanisms include hyperviscosity in the
horizontal, monotone remap in the vertical and dissipation inherent in
the Runge-Kutta method.  HOMME-NH is currently being
evaluated with a suite of test cases from the DCMIP.  
Future work will be to add support
for the height coordinate and coupling with the Energy Exascale Earth
System Model's suite of atmospheric parameterizations.

\appendix
\section{Quasi-Hamiltonian form}
\label{S:appendix}
Following \citeA{DubosTort2014}, we show that 
\eqref{E:uequation}-\eqref{E:continuity} is in quasi-Hamiltonian form plus the
addition of an energetically neutral vertical transport
term $\mathcal{V} \cdot \dot s$.  
Let
\[
\vec{x} =
\begin{pmatrix}
  \uvec \\ w \\ \phi \\ \Theta \\ \dpids 
\end{pmatrix}    
\] and take
\[
\mathcal{H}(\vec{x}) = \iint P + I + K \, dA \, ds
\]
with $P,I,K, dA$ and $ds$ defined as in \S~\ref{S:EnergyConservation}.
To find the functional derivatives of $\mathcal{H}$ we differentiate 
with respect to time as in \S~\ref{S:EnergyConservation} to obtain
\begin{multline}
\frac{\partial \mathcal{H}}{\partial t} = 
\iint
 \dpids \uvec \cdot \frac{\partial \uvec}{\partial t} + 
 \dpids w \frac{\partial w}{\partial t} +  
\left( \dpids - \frac{\partial p}{\partial s} \right) \frac{\partial \phi}{\partial t} + \\
c_p \Pi    \frac{\partial \Theta}{\partial t} + 
 \left( \frac{\uvec^2 + w^2}{2}   +   \phi  \right) \frac{\partial }{\partial t} \dpids  \,dA\,ds 
 \end{multline}
and thus
\begin{align}
\frac{\delta \mathcal{H}}{\delta \uvec} &= \dfrac{\partial \pi}{\partial s}\uvec
\\
\frac{\delta \mathcal{H}}{\delta w} &=\dfrac{\partial \pi}{\partial s}w
\\
\frac{\delta \mathcal{H}}{\delta \phi} &=\dfrac{\partial \pi}{\partial s}-\dfrac{\partial p}{\partial s} = \dfrac{\partial \pi}{\partial s} (1-\mu)\\
\frac{\delta \mathcal{H}}{\delta \Theta} &=c_p\Pi
\\
\frac{\delta \mathcal{H}}{ \delta \tdpids} &=\frac{\uvec^2 + w^2}{2} +\phi.
\end{align}
Using these functional derivatives, the equations of motion can be re-written as
\begin{align*}
&\frac{\partial \uvec}{\partial t} 
+\left(\dfrac{\partial \pi}{\partial s}\right)^{-1}(\nabla_s\times\uvec+f\khat)\times \frac{\delta \mathcal{H}}{\delta \uvec} 
-\left(\dfrac{\partial \pi}{\partial s}\right)^{-1}\frac{\delta \mathcal{H}}{\delta w}\nabla_s w 
-\left(\dfrac{\partial \pi}{\partial s}\right)^{-1}\frac{\delta \mathcal{H}}{\delta \phi}\nabla_s\phi\\
&\quad\quad\quad+\theta_v\nabla_s\frac{\delta \mathcal{H}}{\delta \Theta}
+\nabla_s\frac{\delta \mathcal{H}}{  \delta \tdpids} + \dot{s} \frac{\partial \uvec}{\partial s}=0
\\
&\frac{\partial w}{\partial t} +
\left(\dfrac{\partial \pi}{\partial s}\right)^{-1} (\nabla_sw)\cdot\frac{\delta \mathcal{H}}{\delta \uvec}
+g\left(\dfrac{\partial \pi}{\partial s}\right)^{-1}\frac{\delta \mathcal{H}}{\delta \phi} + \dot{s} \frac{\partial w}{\partial s}=0
\\
&\frac{\partial \phi}{\partial t} 
+\left(\dfrac{\partial \pi}{\partial s}\right)^{-1}(\grad_s\phi)\cdot\frac{\delta \mathcal{H}}{\delta \uvec}
-g\left(\dfrac{\partial \pi}{\partial s}\right)^{-1}\frac{\delta \mathcal{H}}{\delta w} + \dot{s} \frac{\partial \phi}{\partial s}=0
\\
&\frac{\partial \Theta}{\partial t} 
+\nabla\cdot\left(\theta_v\frac{\delta \mathcal{H}}{\delta \uvec}\right) + \frac{\partial}{\partial s} (\dot{s} \Theta)=0
\\
&\frac{\partial}{\partial t}\left(\dpids\right)
+\nabla\cdot\left(\frac{\delta \mathcal{H}}{\delta \uvec}\right) + \frac{\partial}{\partial s} (\dot{s} \dpids)=0
\end{align*}
This can be seen to fall into the quasi-Hamiltonian form
\[
\frac{\partial \vec{x}}{\partial  t} +  \mathcal{J}
\frac{\delta \mathcal{H}}{\delta \vec{x}}
+\mathcal{V} \cdot \dot s = 0
\]
with $\mathcal{J}$ a skew symmetric operator defined via
\begin{equation}
\left( \begin{matrix}
\frac{\mathbf{\eta_s}}{\dpids} \times & -\frac{1}{\dpids} \nabla_s w & -\frac{1}{\dpids} \nabla_s \phi & \theta \nabla_s & \nabla_s \\
\frac{1}{\dpids} \nabla_s w & 0 & \frac{g}{\dpids} & 0 & 0\\
\frac{1}{\dpids} \nabla_s \phi & -\frac{g}{\dpids} & 0 & 0 & 0\\
\nabla_s \cdot \theta & 0 & 0 & 0 & 0\\
\nabla_s \cdot & 0 & 0 & 0 & 0
\end{matrix} \right)
\end{equation}
where $\eta_s = \nabla_s\times\uvec+f\khat$ and 
\begin{equation}
\mathcal{V} \cdot \dot s = \left( \begin{matrix}
\dot s \frac{ \partial \mathbf{u}}{\partial s}\\
\dot s \frac{ \partial w}{\partial s}\\
\dot s \frac{\partial \phi}{\partial s}\\
\frac{\partial}{\partial s} ( \dot s \Theta ) \\
\frac{\partial}{\partial s} (\dot s  \dpids )
\end{matrix} \right)
\end{equation}
We call this quasi-Hamiltonian form since no attempt has been made to prove the Jacobi identity for $\mathcal{J}$. Energy is conserved due the skew-symmetry of $\mathcal{J}$ 
and the condition 
\begin{equation}
\mathcal{V}^* \cdot \frac{\delta \mathcal{H}}{\delta x} = 0
\label{E:VRS}
\end{equation}
where
\begin{equation}
\mathcal{V}^* = \left( \begin{matrix}
\frac{\partial \mathbf{u}}{\partial s}\\
\frac{\partial w}{\partial s}\\
\frac{\partial \phi}{\partial s}\\
- \Theta \frac{\partial}{\partial s}  \\
- \dpids \frac{\partial}{\partial s}
\end{matrix} \right)
\end{equation}
is the adjoint of $\mathcal{V}$. 
The \eqref{E:VRS} condition is a consequence of the vertical relabeling
symmetry of the flow and implies the vertical transport is energetically neutral for
any velocity $\dot s$ with boundary conditions $\dot s =0$.

A general mimetic discretization will preserve the skew-symmetry of $\mathcal{J}$,
but additional properties are needed to obtain a discrete version of 
\eqref{E:VRS}. To see this, we insert the functional derivatives 
into 
\begin{equation}
\label{vert-relabelling}
\frac{\delta \mathcal{H}}{\delta \mathbf{u}} \frac{\partial \mathbf{u}}{\partial s} + \frac{\delta \mathcal{H}}{\delta w} \frac{\partial w}{\partial s} + \frac{\delta \mathcal{H}}{\delta \phi} \frac{\partial \phi}{\partial s} - \Theta \frac{\partial}{\partial s} \frac{\delta \mathcal{H}}{\delta \Theta} - \dpids \frac{\partial}{\partial s} \frac{\delta \mathcal{H}}{\delta \dpids} = 0
\end{equation}
to obtain 
\begin{equation}
\dpids \mathbf{u} \frac{\partial \mathbf{u}}{\partial s} + \dpids w \frac{\partial w}{\partial s} + \dpids (1- \mu) \frac{\partial \phi}{\partial s} - \Theta \frac{\partial}{\partial s} (c_p \Pi) - \dpids \frac{\partial}{\partial s} \left(\frac{\mathbf{u} \cdot \mathbf{u} + w^2}{2} + \phi\right) = 0
\end{equation}
The product rule for $\partial / \partial s$  immediately gives cancellation for the $\mathbf{u}$, $w$ and $\phi$ terms, leaving only
\begin{equation}
-\dpids \mu \frac{\partial \phi}{\partial s} - \Theta \frac{\partial (c_p \Pi)}{\partial s} = 0
\end{equation}
Now recalling that $\dpids \mu = \frac{\partial p}{\partial s}$ and 
\begin{equation}
\frac{1}{\rho} \frac{\partial p}{\partial s} - \theta_v \frac{\partial c_p \Pi}{\partial s} = 0
\end{equation}
which holds in the continuum due to the definition of $\Pi$ and $\theta_v$.  
We can thus obtain a discrete version of \eqref{E:VRS} with the
the SB81 product rule \eqref{E:productrule1} and 
the $\theta_v$ special averaging in \eqref{E:tildedefinition}.